\documentclass[11pt,leqno,twoside]{article}

\usepackage{amsfonts,amsmath,amsthm,amssymb}
\usepackage{enumerate}
\usepackage{graphics,graphicx,subfigure}

\usepackage{color}
\usepackage{todonotes}
\usepackage{cancel}
\usepackage{url}
\usepackage{hyperref}
\usepackage{makeidx}
\usepackage{showidx}
\usepackage{multicol}        
\usepackage{xspace}
\usepackage{stmaryrd}        
\usepackage{pifont}          
\usepackage{fancybox}        
\usepackage{bm}
\usepackage{subfigure}

 \usepackage{fancyhdr}
\theoremstyle{plain}
 \theoremstyle{definition}
 \newtheorem{lem}{Lemma}
 \newtheorem{defn}[lem]{Definition}
 \newtheorem{thm}[lem]{Theorem}
 \newtheorem{prop}[lem]{Proposition}
 \newtheorem{cor}[lem]{Corollary}
 \newtheorem{notn}[lem]{Notations}
 \newtheorem{pb}[lem]{Problem}
 \newtheorem{form}[lem]{Formulae}
 
 \newtheorem*{rk}{Remark}
 \newtheorem*{com}{Comment}
 \newtheorem*{ex}{Example}
 \theoremstyle{remark}

 \newcommand{\blem}{\begin{lem}}
 \newcommand{\elem}{\end{lem}}
 \newcommand{\bdefn}{\begin{defn}}
 \newcommand{\edefn}{\end{defn}}
 \newcommand{\bthm}{\begin{thm} }
 \newcommand{\ethm}{\end{thm}}
 \newcommand{\bprop}{\begin{prop}}
 \newcommand{\eprop}{\end{prop}}
 \newcommand{\bcor}{\begin{cor}}
 \newcommand{\ecor}{\end{cor}}
 \newcommand{\bnotn}{\begin{notn}}
 \newcommand{\enotn}{\end{notn}}
 \newcommand{\bpb}{\begin{pb}}
 \newcommand{\epb}{\end{pb}}
 \newcommand{\bform}{\begin{form}}
 \newcommand{\eform}{\end{form}}
 \newcommand{\brk}{\begin{rk}}
 \newcommand{\erk}{\end{rk}}
 \newcommand{\bcom}{\begin{com}}
 \newcommand{\ecom}{\end{com}}
 \newcommand{\bex}{\begin{ex}}
 \newcommand{\eex}{\end{ex}}
 \newcommand{\bpf}{\begin{proof}}
 \newcommand{\epf}{\end{proof}}

\newcommand{\vp}{{\bf p}}

\newcommand{\vy}{{\bf y}}

\newcommand{\vJ}{{\bf J}}

\newcommand{\vQ}{{\bf Q}}
\newcommand{\vR}{{\bf R}}

\newcommand{\vW}{{\bf W}}

\newcommand{\cP}{\mathcal{P}}

\newcommand{\cT}{\mathcal{T}}
\newcommand{\cU}{\mathcal{U}}
\newcommand{\cV}{\mathcal{V}}

\newcommand{\bC}{\mathbb{C}}

\newcommand{\bN}{\mathbb{N}}

\newcommand{\bR}{\mathbb{R}}

\newcommand{\be}{\begin{equation}}
\newcommand{\ee}{\end{equation}}
\newcommand{\bal}{\begin{align}}
\newcommand{\eal}{\end{align}}
\newcommand{\ba}{\begin{align*}}
\newcommand{\ea}{\end{align*}}
\newcommand{\bmx}{\begin{matrix}}
\newcommand{\emx}{\end{matrix}}
\newcommand{\bbmx}{\begin{bmatrix}}
\newcommand{\ebmx}{\end{bmatrix}}
\newcommand{\bpmx}{\begin{pmatrix}}
\newcommand{\epmx}{\end{pmatrix}}
\newcommand{\bvmx}{\begin{vmatrix}}
\newcommand{\evmx}{\end{vmatrix}}

\newcommand{\wt}{\widetilde}
\newcommand{\f}{\frac}
\newcommand{\df}{\dfrac}

\newcommand{\imp}{\Longrightarrow}

\newcommand{\inc}{\subseteq}

\newcommand{\tto}{\longrightarrow}

\newcommand{\sgn}{\mathrm{sgn}}

\newcommand{\argmin}{{\rm argmin}\,}

\newcommand{\minimize}[1]{\underset{#1}{\rm minimize}\,}

\newcommand{\eps}{\varepsilon}
  
\newcommand{\Toep}{{\rm \bf Toep}} 
\newcommand{\rev}[1]{#1}

\title{\vspace{-25mm}Computation of  Chebyshev Polynomials \rev{for} Union of Intervals\medskip\hrule height 1.2pt \vspace{-6mm}}
\author{Simon Foucart\footnote{The research of the first author is partially funded by the NSF grants DMS-1622134 and DMS-1664803.} \,(Texas A\&M University) and Jean Bernard Lasserre\footnote{The research of the second author is funded by the European Research Council (ERC) under
the \rev{European Union's} Horizon 2020 research and innovation program (grant agreement 666981 TAMING).} \,(LAAS-CNRS)}
\date{\vspace{-6mm}\rule{100mm}{0.8pt}}

\newcommand\shorttitle{Computation of  Chebyshev Polynomials on Union of Intervals}
\newcommand\authors{S. Foucart, J. B. Lasserre}

\begin{document}
\maketitle

\vspace{-15mm}
\begin{abstract}
Chebyshev polynomials of the first and second kind \rev{for} a set $K$ are monic polynomials with minimal \rev{$L_\infty$- and $L_1$-norm} on $K$,
respectively.
This articles presents numerical procedures based on semidefinite programming to compute these polynomials in case $K$ is a finite union of compact intervals.
For Chebyshev polynomials of the first kind,
the procedure makes use of a characterization of polynomial nonnegativity.
It can incorporate additional constraints,
e.g. that all the roots of the polynomial lie \rev{in} $K$.
For Chebyshev polynomials of the second kind,
the procedure exploits the method of moments.
\end{abstract}



\vspace{-5mm}
\begin{center}
\rule{100mm}{0.8pt}
\end{center}


\section{Introduction}

The $N$th Chebyshev polynomial \rev{for} a compact \rev{infinite} subset $K$ of $\bC$  is defined as the monic polynomial of degree $N$ with minimal max-norm on $K$.
Its uniqueness is a straightforward consequence of the uniqueness of best polynomial approximants to a continuous function (here $z\mapsto z^N$) with respect to the max-norm, see e.g. \cite[p. 72, Theorem 4.2]{CA}.
We shall denote it as $\cT_N^K$, i.e.,
\be
\cT_N^{K}
= \underset{P(z) = z^N + \cdots }{\argmin} \|P\|_K,
\qquad \mbox{where }
\|P\|_K = \max_{z \in K} |P(z)|.
\ee
We reserve the notation $T_N^K$ for the Chebyshev polynomial normalized to have max-norm equal to one on $K$, i.e.,
\be
T_N^K = \f{\cT_N^K}{ \| \cT_N^K  \|_K}.
\ee
With this notation, the usual $N$th Chebyshev polynomial (of the first kind) satisfies
\be
T_N = T_N^{[-1,1]} = 2^{N-1} \cT_N^{[-1,1]},
\qquad N \ge 1. 
\ee
Chebyshev polynomials \rev{for} a compact subset $K$ of $\bC$ play an important role in \rev{logarithmic} potential theory.
For instance, it is known that the capacity ${\rm cap}(K)$ of $K$
is related to the Chebyshev numbers $t_N^K := \|\cT_N^K\|_K$ via 
\be
\label{Conv}
\left( t_N^K \right)^{1/N} \underset{N \to \infty}{\tto} {\rm cap}(K),
\ee
see \cite[p.163, Theorem 3.1]{ST} for a weighted version of this statement.
The articles \cite{CSZ,CSYZ} recently studied in greater \rev{detail} the asymptotics of the convergence \eqref{Conv} in case $K$ is a subset of~$\bR$.
\rev{This being said}, the capacity is in general hard to determine
--- it can be found explicitly in a few specific situations,
e.g. when $K$ is the inverse image of an interval by certain polynomials
(see \cite[Theorem~11]{GvA}),
and otherwise some numerical \rev{methods} for computing the capacity have been proposed in \cite{RR},
\rev{see also Section 5.2 of \cite{Ran}}.
As for the Chebyshev polynomials, one is tempted to anticipate a worse state of affairs.
However, this is not the case for the situation considered in this article,
i.e., when $K \subseteq [-1,1]$ is a finite union of $L$ compact intervals\footnote{\rev{The assumption $K \subseteq [-1,1]$ is not restrictive, 
as any compact subset of $\bR$ can be moved into the interval $[-1,1]$ by an affine transformation.}},
say 
\be
\label{K}
K = \bigcup_{\ell = 1}^L [a_\ell ,b_\ell],
\qquad \quad -1=a_1 < b_1 < a_2 < b_2 < \cdots < a_L < b_L=1.
\ee
There are explicit constructions of Chebyshev polynomials (as orthogonal polynomials with a predetermined weight, see \cite[Theorem 2.3]{Peh}),
albeit only under the condition that $\cT_N^K$ is a {\em strict} Chebyshev polynomial
(meaning that it possesses $N+L$ points of equioscillation on $K$
--- a condition which is verifiable a priori, see \cite[Theorem 2.5]{Peh}).
Chebyshev polynomials can otherwise be computed using Remez-type algorithms \rev{for} finite unions of intervals, see  \cite{Fil}.

\rev{A first contribution} of this article is to put forward an alternative numerical procedure  that enables the \rev{accurate} computation of the Chebyshev polynomials
whenever $K$ is a finite union of compact intervals.
The procedure, based on semidefinite programming as described in Section \ref{Sec1stKind}, can also incorporate a \rev{weight $w$ (i.e., a continuous and positive function on $K$), restricted here to be a rational function,}
and output the polynomials
\be
\label{DefTw}
\cT_N^{K,w} 
= \underset{P(x) = x^N + \cdots }{\argmin} \; \bigg\| \f{P}{w}  \bigg\|_K.
\ee
An appealing feature of this approach is that extra constraints can easily be incorporated in the minimization of \eqref{DefTw}.
For instance, we will show how to compute the $N$th {\em restricted} Chebyshev polynomial on $K$,
i.e., the monic polynomial of degree $N$ having all its roots in $K$ with minimal max-norm on $K$.

\rev{A second contribution} of this article is to propose another semidefinite-programming-based procedure to compute weighted Chebyshev polynomial\rev{s} of the second kind, so to \rev{speak}.
By this, we mean  polynomials\footnote{\rev{The uniqueness of $\cU_N^{K,w}$ is not necessarily guaranteed: in the unweighted case, one can e.g. check that the monic linear polynomials with minimal $L_1$-norm on $K=[-1,-c]\cup [c,1]$ are all the $x-d$, $d \in [-c,c]$.
We will not delve into conditions ensuring uniqueness of $\cU_N^{K,w}$ in this article.}}
\be
\cU_N^{K,w} 
\rev{\in} \underset{P(x) \in x^N + \cdots }{\argmin} \; \bigg\| \f{P}{w}  \bigg\|_{L_1(K)}.
\ee
\rev{The restriction that the weight $w$ is a rational function is not needed here},
but this time the computation is only approximate.
Nonetheless, it produces lower and upper bounds for the genuine minimium $\|\cU_N^{K,w}/w\|_{L_1(K)}$.
Both bounds are proved to converge to the genuine minimum as a parameter $d \ge N$ grows to infinity.
\rev{Along the way},
we shall prove that the Chebyshev polynomial of the second kind \rev{for} $K$\rev{, if unique}, has \rev{simple roots all} lying inside $K$.

The procedures for computing Chebyshev polynomials of the first and second kind have been implemented in {\sc matlab}.
They rely on the external packages {\sf CVX} (for specifying and solving convex programs \cite{CVX}) 
and {\sf Chebfun} (for numerically computing with functions \cite{Chebfun}).
They can be downloaded from the authors' webpage as part of the reproducible file accompanying this article.

\section{Chebyshev polynomials of the first kind}
\label{Sec1stKind}

With $K$ as in \eqref{K},
we consider a rational\footnote{We could also work with piecewise rational weight functions, but we choose not to \rev{do} so in order to avoid overloading already heavy notation.} weight function $w$ taking the form
\be
\label{w}
w = \f{\Sigma}{\Omega},
\ee
where the polynomials $\Sigma$ and $\Omega$ are positive on each $[a_\ell,b_\ell]$. 
We shall represent polynomials $P$ of degree at most $N$ by their Chebyshev expansions written as
\be
P = \sum_{n=0}^N p_n T_n.
\ee
In this way, finding the $N$th Chebyshev polynomial of the first kind \rev{for} $K$ with weight $w$
amounts to solving the optimization problem
\be
\minimize{p_0, p_1,\ldots, p_N \in \bR}  \max_{\ell=1:L} \left\| \f{\Omega P}{\Sigma}  \right\|_{[a_\ell,b_\ell]} 
\qquad \mbox{s.to} \quad p_N = \f{1}{2^{N-1}}.
\ee
After introducing a slack variable $c \in \bR$,
this is equivalent to the optimization problem
\be
\minimize{c, p_0, p_1,\ldots, p_N \in \bR} \; \;  c
\qquad \mbox{s.to}\quad
p_N = \f{1}{2^{N-1}}
\quad \mbox{and} \quad
\left\| \f{\Omega P}{\Sigma}  \right\|_{[a_\ell,b_\ell]} \le c 
\quad \mbox{for all }\ell = 1:L.
\ee
The latter constraints can be rewritten as 
$-c \le \Omega P / \Sigma \le c$ on $[a_\ell,b_\ell ]$, $\ell=1:L$,
i.e., as the two polynomial nonnegativity constraints
\be
c \Sigma(x) \pm \Omega(x) P(x) \ge 0
\qquad \quad \mbox{for all } x \in [a_\ell, b_\ell ]  \quad \mbox{ and all} \quad \ell = 1:L.
\ee
The key to the argument is now to exploit an exact semidefinite characterization of these constraints.
This is based on the following result,
 which was \rev{established and utilized} in \cite{FP}, see Theorem~3 there.

\bprop
Given $[a,b] \inc [-1,1]$ and a polynomial $C(x) = \sum_{m=0}^M c_m T_m(x)$ of degree at most~$M$,
the nonnegativity condition
\be
\label{Pos}
C(x) \ge 0 
\qquad \mbox{for all }x \in [a,b]
\ee
is equivalent to the existence of semidefinite matrices $\vQ \in \bC^{(M+1)\times (M+1)}$, $\vR \in \bC^{M \times M}$ such that
\be 
\label{PosChar}
 \sum_{i-j = m} Q_{i,j} 
+ \alpha  \sum_{i-j=m-1} R_{i,j} 
-\beta \sum_{i-j=m} R_{i,j} 
+ \overline{\alpha}  \sum_{i-j=m+1} R_{i,j} 
\; = \; \left\{ \bmx \f{1}{2}c_m, & m = 1:M \\
c_0, & m=0 \emx \right\},
\ee
where $ \alpha = \f{1}{2}\exp\left( \f{\imath}{2} \arccos(a)+ \f{\imath}{2}\arccos(b) \right)$ and $\beta = \cos \left( \f{1}{2}\arccos(a) - \f{1}{2}\arccos(b) \right)$.
\eprop

In the present situation,
we apply this result to the polynomials $C = c \Sigma \pm \Omega P$ required to be nonnegative on each $[a_\ell,b_\ell]$.
With 
\be
M:= \max \left\{ \deg(\Sigma), \deg(\Omega) + N \right\},
\ee
we write the Chebyshev expansions of $\Sigma$ and of $\Omega P$ as
\be
\Sigma = \sum_{m=0}^{M} \sigma_m T_m,
\qquad \qquad
\Omega P = \sum_{m=0}^{M} (\vW \vp)_m T_m,
\ee
where $\vW \in \bR^{(M+1) \times (N+1)}$ is the matrix of the linear map transforming the Chebyshev coefficients of $P$ into the Chebyshev coefficients of $\Omega P$.
Our considerations can now be summarized as follows.

\bthm
The $N$th Chebyshev polynomial $\cT_N^{K,w}$ \rev{for} the set $K$ given in \eqref{K} and with weight $w$ given in \eqref{w} has Chebyshev coefficients $p_0,p_1,\ldots,p_N$ that solve the semidefinite program
\begin{align}
\label{SDP1Kd}
\minimize{\substack{c,p_0,p_1,\ldots,p_N \in \bR \\ 
\vQ^{\pm, \ell} \in \bR^{(M+1)\times (M+1)}\\
\vR^{\pm, \ell} \in \bR^{M \times M} } }
\; c
& & &\mbox{s.to}
\quad p_N=\f{1}{2^{N-1}}, 
\quad \vQ^{\pm, \ell} \succeq {\bf 0},
\quad \vR^{\pm, \ell} \succeq {\bf 0},\\
\nonumber
& & & \mbox{and} \quad 
\sum_{i-j = m} Q^{\pm,\ell}_{i,j} +  \alpha_\ell \sum_{i-j=m-1} R^{\pm,\ell}_{i,j} 
-\beta_\ell \sum_{i-j=m} R^{\pm,\ell}_{i,j}
+ \overline{\alpha_\ell} \sum_{i-j=m+1} R^{\pm,\ell}_{i,j}\\
\nonumber
& & & \qquad \qquad \qquad \quad \; \,
=\left\{ \bmx \f{1}{2} \sigma_m c \pm \f{1}{2} (\vW \vp)_m , & m = 1:M \\
\sigma_0 c \pm (\vW \vp)_0, & m=0 \emx \right\},
\end{align}
where
$ \alpha_\ell = \f{1}{2}\exp\left( \f{\imath}{2} \arccos(a_\ell)+ \f{\imath}{2}\arccos(b_\ell) \right)$ and $\beta_\ell = \cos \left( \f{1}{2}\arccos(a_\ell) - \f{1}{2}\arccos(b_\ell) \right)$.
\ethm

\begin{figure}[htbp]
\center
\subfigure[\rev{$K=K_1$}
]{
\includegraphics[width=0.48\textwidth,height=0.33\textwidth]{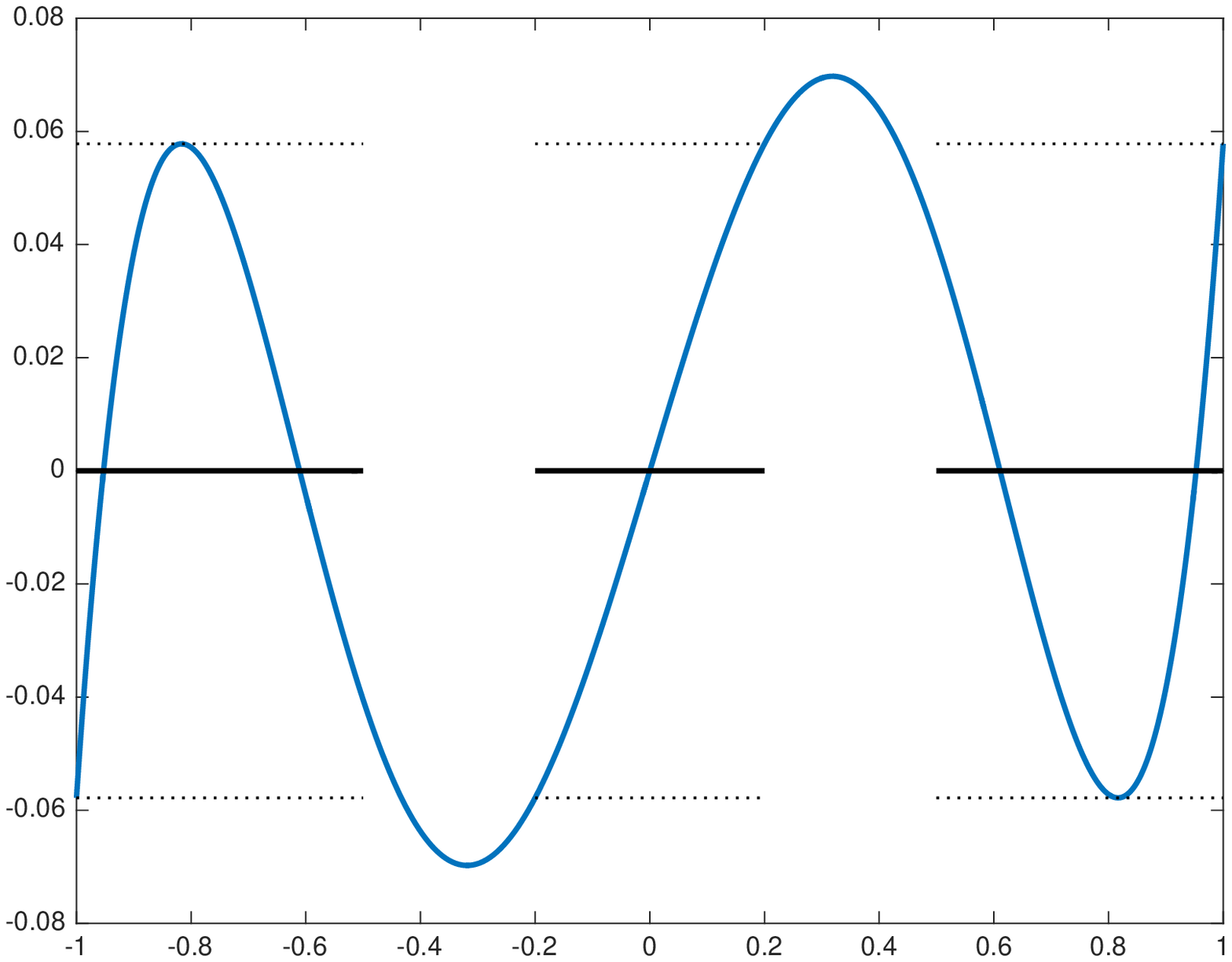}
}
\subfigure[
\rev{$K=K_1$, weighted}
]{
\includegraphics[width=0.48\textwidth,height=0.33\textwidth]{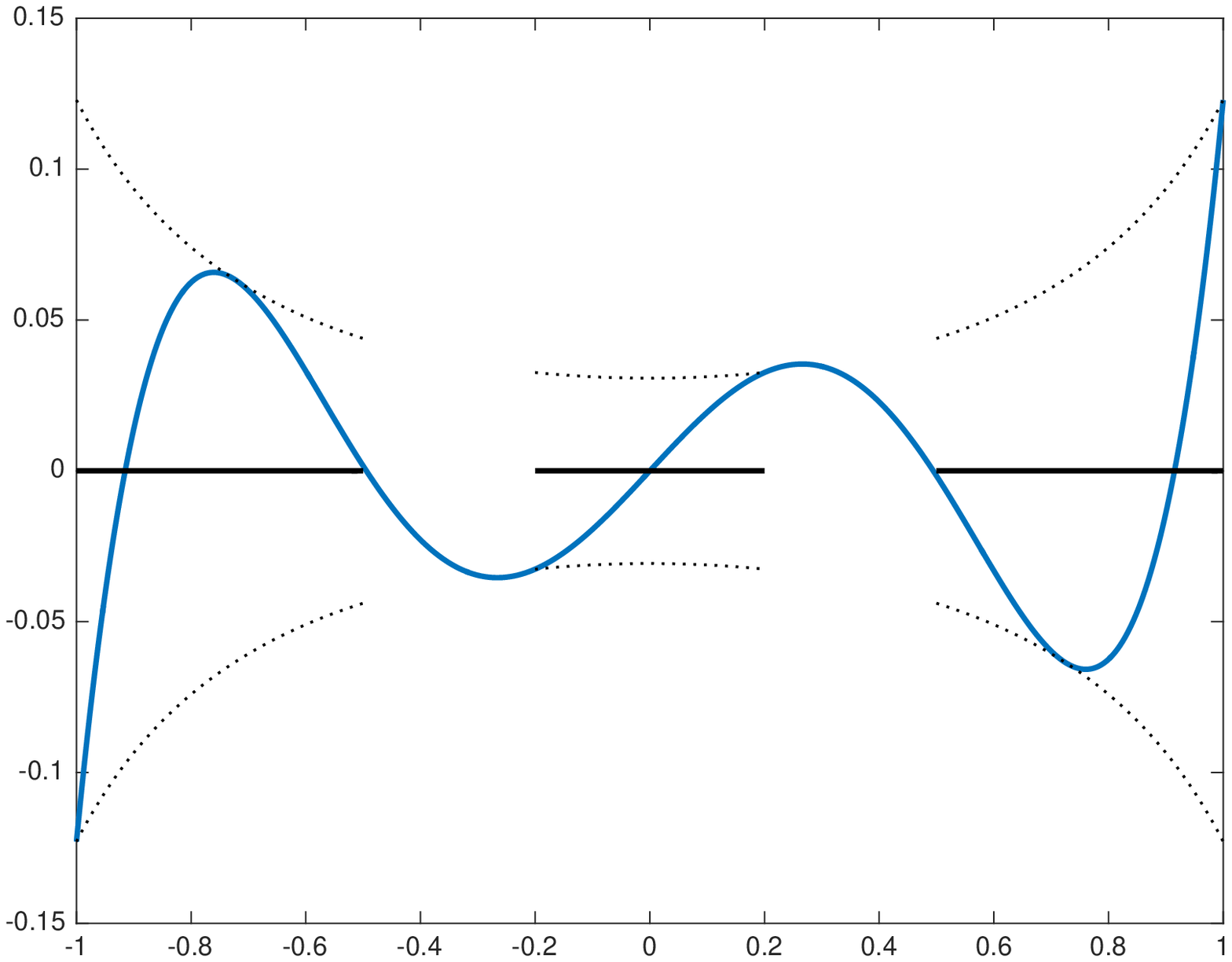}
}

\subfigure[
\rev{$K=K_2$}
]{
\includegraphics[width=0.48\textwidth,height=0.33\textwidth]{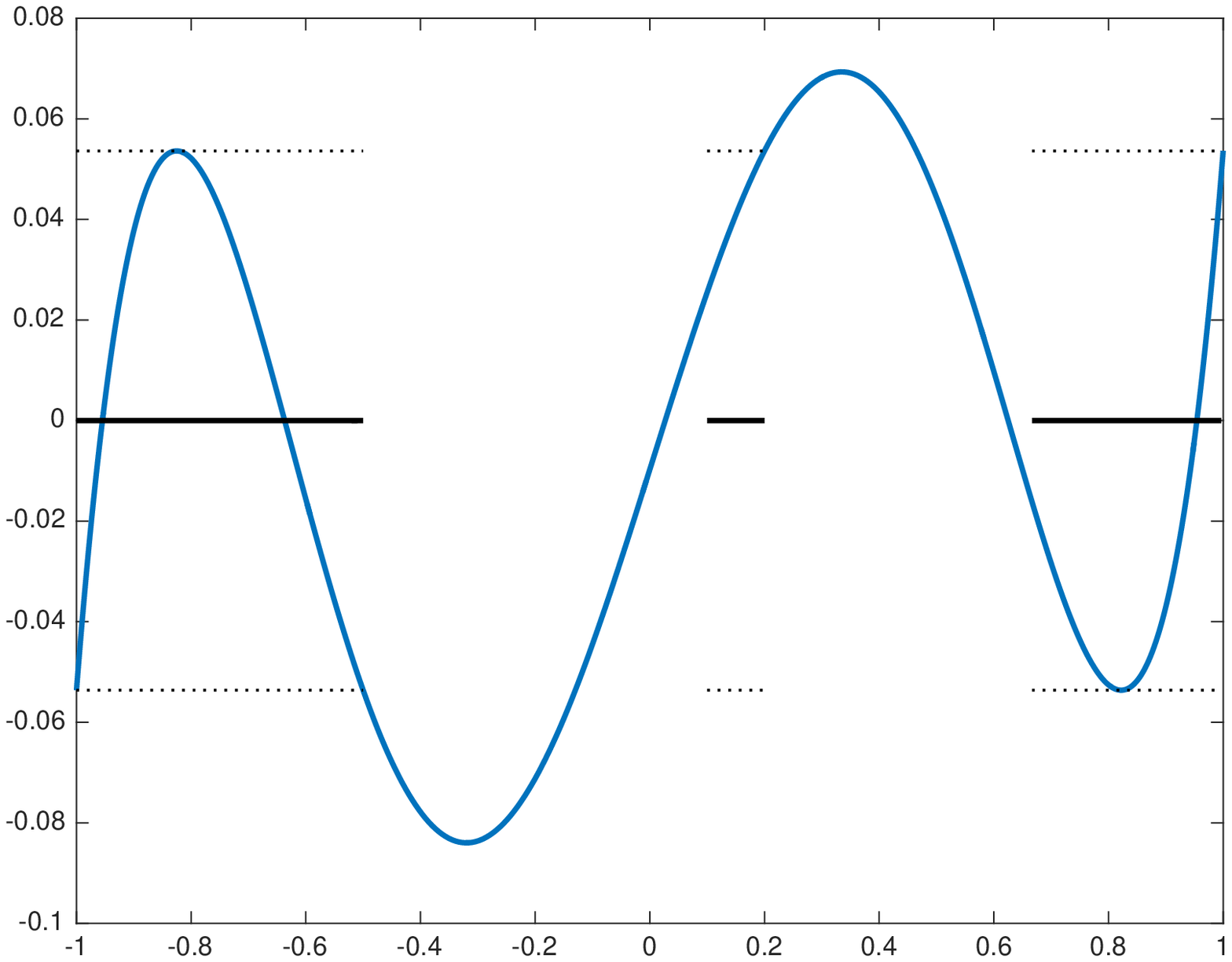}
}
\subfigure[
\rev{$K=K_2$, weighted}
]{
\includegraphics[width=0.48\textwidth,height=0.33\textwidth]{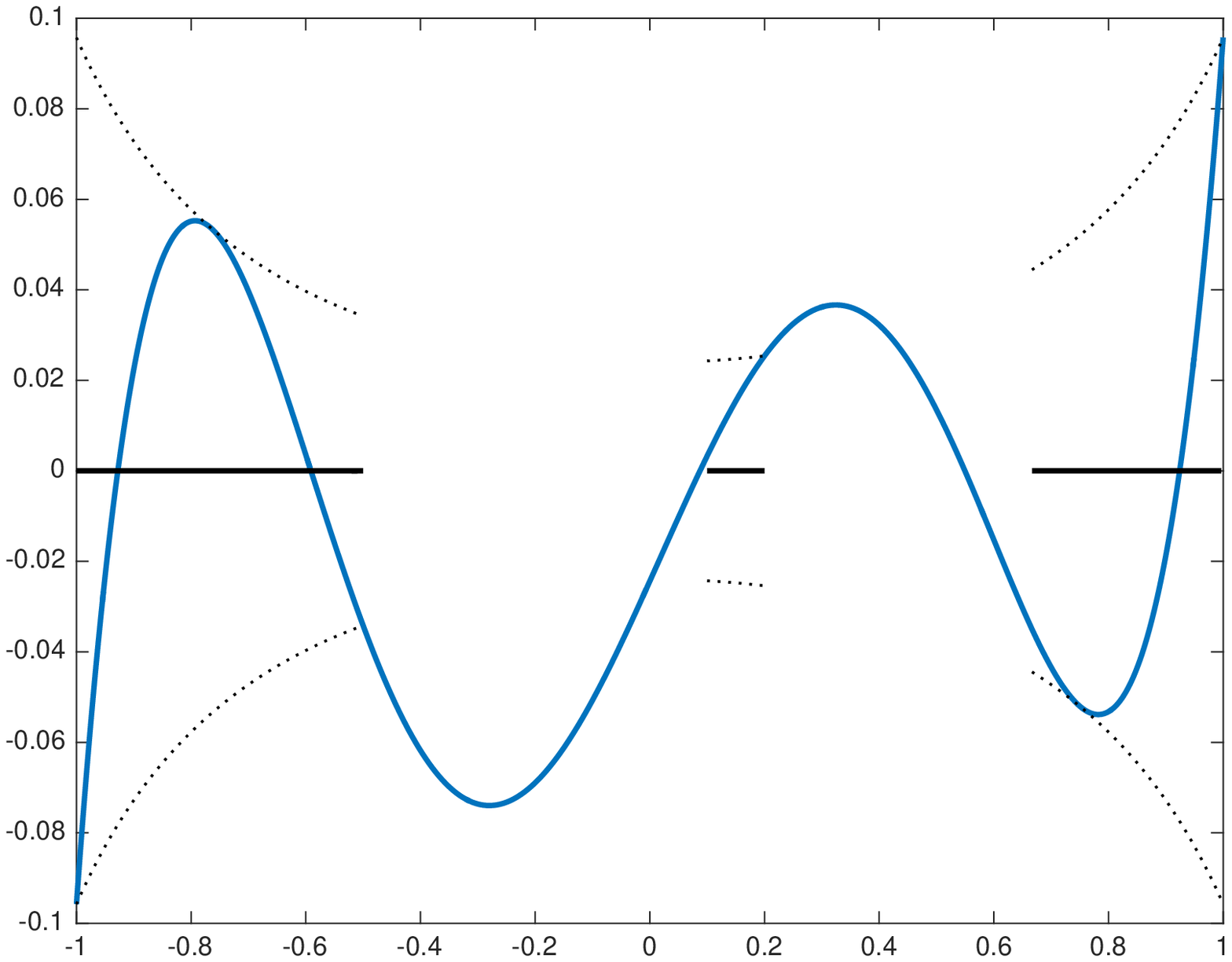}
}

\subfigure[
\rev{$K=K_2$, restricted}
]{
\includegraphics[width=0.48\textwidth,height=0.33\textwidth]{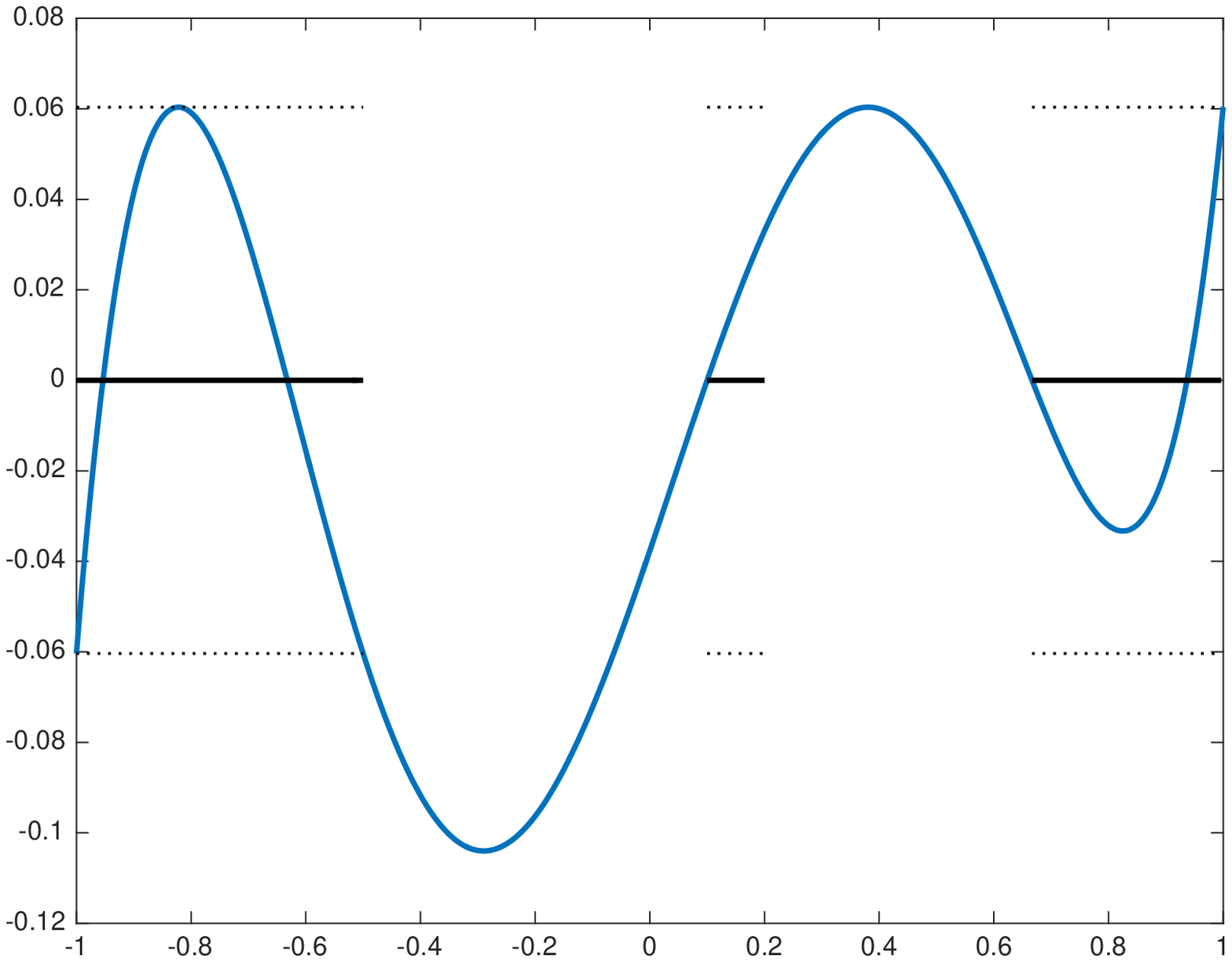}
}
\subfigure[
\rev{$K=K_2$, restricted and weighted}
]{
\includegraphics[width=0.48\textwidth,height=0.33\textwidth]{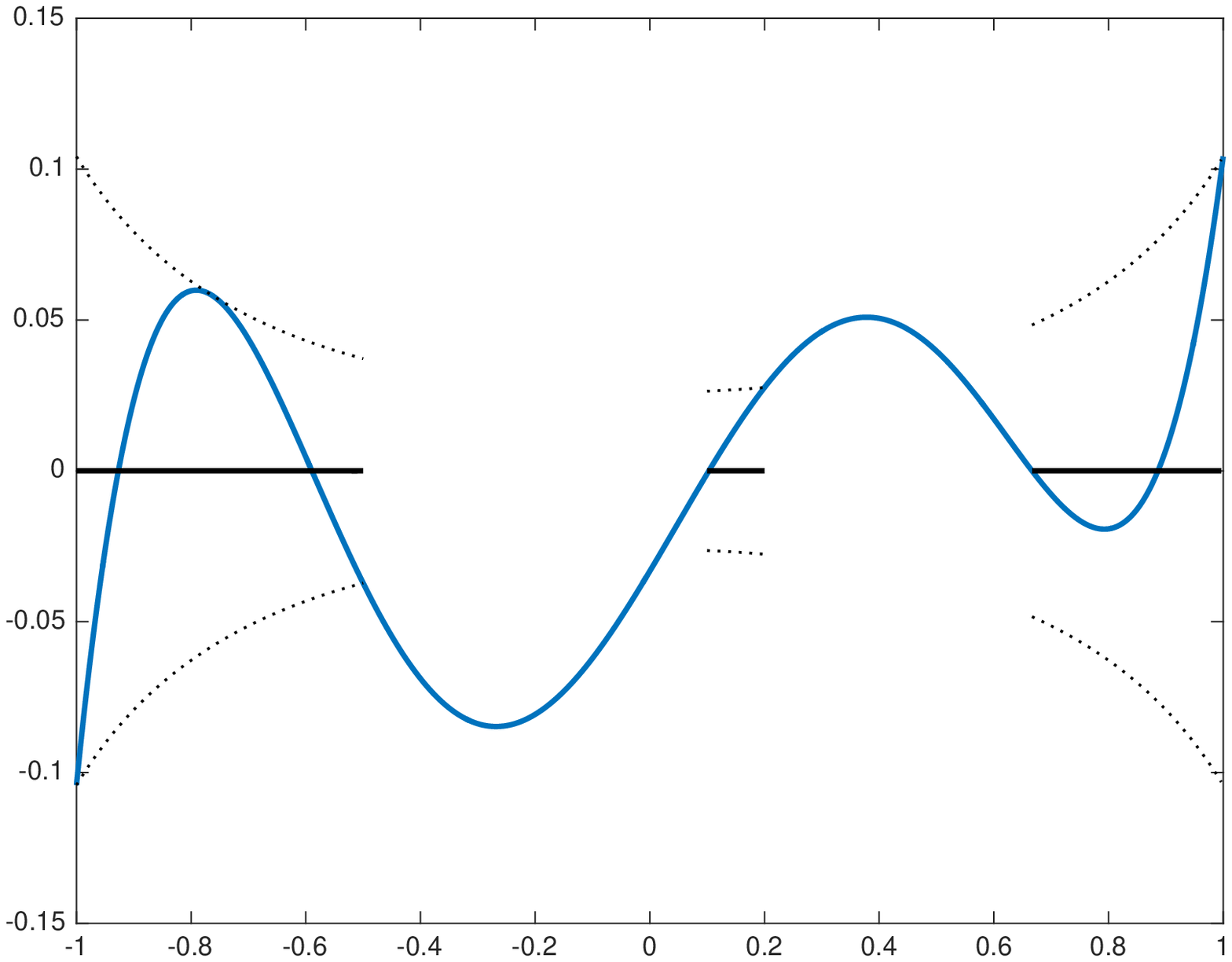}
}
\caption{\rev{$5$th Chebyshev polynomials of the first kind for $K_1 = \big[-1,-\f{1}{2}\big] \cup \big[-\f{1}{5},\f{1}{5} \big] \cup \big[\f{1}{2},1 \big]$ and for $K_2 = \big[ -1,-\f{1}{2} \big] \cup \big[\f{1}{10},\f{1}{5}\big] \cup \big[\f{2}{3},1\big]$:
the first two rows correspond to the unrestricted case,
while restricted Chebyshev polynomials are shown in the last row;
the first column corresponds to the unweighted case,
while weighted Chebyshev polynomials with weight $w(x) = (1+x^2)/(2-x^2)$ are shown in the second column.}
}
\label{Fig1}
\end{figure}

Figure \ref{Fig1} provides examples of Chebyshev polynomials of degree $N=5$ \rev{for} the union of $L=3$ intervals which were computed by solving \eqref{SDP1Kd}.
In all cases, the Chebyshev polynomials equioscillate $N+1~=~6$ times between $-w$ and $+w$ on $K$, as they should.
However, they are not strict Chebyshev polynomials,
since the number of equioscillation points on $K$ is smaller than $N+L=8$.
We notice in (c) and (d) that some roots of the Chebyshev polynomials do not lie in the  set $K$.
We display in (e) and (f) the {\em restricted} Chebyshev polynomial \rev{for} $K$,
i.e., the monic polynomial of degree $N$ with minimal max-norm on $K$ which satisfies the additional constraint that all its roots lie in $K$.
This constraint reads 
\be
\label{ResCons}
P
\mbox{ does not vanish on }(b_\ell, a_{\ell+1}),
\quad \ell = 1:L-1.
\ee
We consider the semidefinite program \eqref{SDP1Kd} supplemented with the relaxed constraint 
\be
\label{ResConsRel}
P
\mbox{ does not change sign on }[b_\ell, a_{\ell+1}],
\quad \ell = 1:L-1.
\ee
This is solved by selecting the smallest value (along with the corresponding minimizer) among the minima of $2^{L-1}$ semidefinite programs \eqref{SDP1Kd} indexed by $(\eps_1,\ldots,\eps_{L-1}) \in \{ \pm 1\}^{L-1}$,
where the added constraint is the semidefinite characterization of the polynomial nonnegativity condition
\be
\eps_\ell P(x) \ge 0
\qquad \mbox{for all } x \in [b_\ell,a_{\ell+1}]
\quad \mbox{and all}\quad \ell = 1:L-1.
\ee
One checks whether the selected minimizer satisfies the original constraint \eqref{ResCons}.
If it does, 
then the restricted Chebyshev polynomial has indeed been found,
as in (e) and (f) of Figure \ref{Fig1}.

\brk
Concerning the computation of the capacity of a union of intervals,
we do not recommend using our semidefinite procedure or a Remez-type procedure to produce Chebyshev polynomials before invoking \eqref{Conv} to approximate the capacity.
If one really wants to take such a route, 
it seems wiser to work with the numerically-friendlier orthogonal polynomials
\be
\cP_N^K  = \underset{P(x) = x^N + \cdots }{\argmin} \; \| P  \|_{L_2(K)}.
\ee
Indeed, we also have
\be
\|\cP_N^K\|_{L_2(K)}^{1/N} \underset{N \to \infty}{\tto} {\rm cap}(K),
\ee
as a consequence of the inequalities
\be
\f{1}{N+1} \bigg[ \min_{\ell = 1:L} (b_\ell - a_\ell) \bigg] ^{1/2}  \|\cT_N^K\|_K 
\le \|\cP_N^K\|_{L_2(K)}
\le \bigg[ \sum\nolimits_{\ell=1}^L (b_\ell - a_\ell) \bigg]^{1/2} \|\cT_N^K\|_K.
\ee
\erk

\section{Chebyshev polynomials of the second kind}
\label{Sec2ndKind}

Still with $K$ as in \eqref{K},
but with an arbitrary \rev{(positive and continuous)} weight function $w$,
we are \rev{now} targeting $N$th Chebyshev polynomial\rev{s} of the second kind \rev{for} $K$ with weight $w$,
i.e.,
\be 
\label{DefU}
\cU_N^{K,w} \rev{\in}  \underset{P(x) = x^N + \cdots }{\argmin} \left\| \f{P}{w} \right\|_{L_1(K)},
\qquad \mbox{where} \quad
\left\| \f{P}{w} \right\|_{L_1(K)}
= \sum_{\ell=1}^L \int_{a_\ell}^{b_\ell} \f{|P(t)|}{w(t)} dt.
\ee
Let us drop the superscript $w$ and simply write $\cU_N^K$ for $\cU_N^{K,w}$.
Minimizing the $L_1$-norm on $K$ exactly seems out of reach,
so instead we shall perform the minimization of a more tractable ersatz norm,
which will be formally defined in Proposition~\ref{ThmErs}.
This ersatz norm stems from a reformulation of the $L_1$-norm on~$K$, as described in the steps below.
Given a polynomial $P$ of degree at most $N$,
we start by making two changes of variables to write
\be
\left\| \f{P}{w} \right\|_{L_1(K)}
= \sum_{\ell = 1}^L \f{b_\ell - a_\ell}{2} \int_{-1}^1 \f{|P_{\ell} (x)|}{w_\ell(x)} dx
= \sum_{\ell = 1}^L \f{b_\ell - a_\ell}{2} \int_{0}^{\pi} |P_{\ell} (\cos(\theta))| \f{\sin(\theta) d\theta}{w_\ell(\cos(\theta))} ,
\ee
where $P_\ell$ and $w_\ell$ denote the functions $P_{|[a_\ell,b_\ell]}$ and $w_{|[a_\ell,b_\ell]}$ transplanted to $[-1,1]$, for instance
\be
P_\ell (x) = P \left(
\f{(b_\ell - a_\ell)x + a_\ell + b_\ell}{2}
\right),
\qquad x \in [-1,1].
\ee
We continue by decomposing the signed measures $P_\ell (\cos(\theta)) \sin(\theta) / w_\ell(\cos(\theta)) d\theta$ as differences of two nonnegative measures, so that 
\be
\label{L1Reform1}
\left\| \f{P}{w} \right\|_{L_1(K)} = 
\inf_{\mu_1^{\pm},\ldots, \mu_L^\pm}
\sum_{\ell = 1}^L \f{b_\ell - a_\ell}{2}
\int_0^\pi d(\mu_\ell^+ + \mu_\ell^-)
\quad \mbox{s.to }
d(\mu_\ell^+ - \mu_\ell^-)(\theta) = P_\ell(\cos(\theta)) \f{\sin(\theta) d\theta}{w_\ell(\cos(\theta))} ,
\ee
where the infimum is taken over all nonnegative measures on $[0,\pi]$.
As is well known, a minimization over nonnegative measures can be reformulated as a minimization over their sequences of moments.
There are several options to do so: here, emulating an approach already exploited in \cite{FL}, see Section 3 there,
we rely on the discrete trigonometric moment problem encapsulated in the following statement.

\bprop
Given a sequence $\vy \in \bR^\bN$,
there exists a nonnegative measure $\mu$ on $[0,\pi]$ such that
\be 
\int_{0}^\pi \cos( k \theta) d \mu(\theta)
 = y_k,
 \qquad 
 k \ge 0,
 \ee
 if and only if the infinite Toeplitz matix build from $\vy$ is positive semidefinite, i.e.,
 \be 
 \Toep_\infty(\vy) := \bbmx
y_0 & y_1 & y_2 & y_3 & \cdots\\
y_1 & y_0 & y_1 & y_2 & \\
\rev{y_2} & y_1 & y_0 & y_1 & \ddots \\
y_3 & y_2 & y_1 & \ddots & \ddots \\
\vdots & & \ddots & \ddots & \ddots
\ebmx
 \succeq 0.
 \ee
 The latter means that all the finite sections of $ \Toep_\infty(\vy)$ are positive semidefinite,
 i.e.,
\be
 \Toep_d(\vy) = 
\bbmx 
y_0 & y_1 & \cdots & \cdots & y_{d}\\
y_{1} & y_{0} & y_{1} & & \vdots \\
\vdots & y_{1} & \ddots & \ddots & \vdots\\
\vdots & & \ddots & \ddots & y_{1}\\
y_{d} & \cdots & \cdots & y_{1} & y_{0}
 \ebmx 
 \succeq 0
 \qquad \mbox{for all } d\ge 0.
\ee
  \eprop

With $\vy^{1,\pm},\ldots,\vy^{L,\pm} \in \bR^\bN$ representing the sequences of moments of $\mu^\pm_1,\ldots,\mu^\pm_L$,
the objective function in \eqref{L1Reform1} just reads \vspace{-1mm}
$$
\sum_{\ell = 1}^L \f{b_\ell - a_\ell}{2}
\int_0^\pi d(\mu_\ell^+ + \mu_\ell^-)
=
\sum_{\ell=1}^L \f{b_\ell - a_\ell}{2} \left( y^{\ell,+}_0 + y^{\ell,-}_0 \right).
$$
As for the constraints in \eqref{L1Reform1},
with $\vW^\ell \in \bR^{(N+1) \times (N+1)}$ denoting the matrix of the linear map transforming the Chebyshev coefficients of $P$ into the Chebyshev coefficients of the second kind of~$P_\ell$,
so that \vspace{-2mm}
\be
P_\ell
= \sum_{n=0}^N (\vW^\ell \vp)_n U_n,
\ee
they become, for all $\ell=1:L$ and all $k \ge 0$, \vspace{-1mm}
\be
y^{\ell,+}_k - y^{\ell,-}_k 
= \sum_{n=0}^N  (\vW^\ell \vp)_n  \int_0^\pi \cos(k \theta) U_n(\cos(\theta)) \f{\sin(\theta) d\theta}{w_\ell(\cos(\theta))}
 = (\vJ^\ell \vW^\ell \vp)_k,
\ee
where the infinite matrices $\vJ^\ell \in \bR^{\bN \times (N+1)}$ have entries \vspace{-1mm}
\be
\label{EntriesJ}
J^\ell_{k,n} =  \int_0^\pi \f{\cos(k \theta) \sin((n+1)\theta)}{w_\ell(\cos(\theta))} d\theta .
\ee 
The finite matrices $\vJ^{\ell,d} \in \bR^{(d+1)\times (N+1)}$, obtained by keeping the first $d+1$ rows of $\vJ^\ell$,
are to be precomputed \rev{numerically} and can  sometimes \rev{even} be determined explicitly, 
e.g. \vspace{-1mm}
\be
\mbox{when $w=1$,} \quad
J^{\ell,d}_{k,n} 
 = \left\{
\bmx
0 & \mbox{if $k$ and $n$ have different parities},\\
\df{2(n+1)}{(n+1)^2-k^2} & \mbox{if $k$ and $n$ have similar parities}.
\emx
\right.
\ee
Taking into account the constraints that the $\vy^{\ell,\pm} \in \bR^\bN$ must be sequences of moments,
we arrive at a semidefinite reformulation of the weighted $L_1$-norm on $K$ given by \vspace{-1mm}
\begin{align}
\label{L1Reform2}
\left\| \f{P}{w} \right\|_{L_1(K)} & =
\inf_{\vy^{1,\pm},\ldots,\vy^{L,\pm} \in \bR^{\bN}}
\sum_{\ell=1}^L \f{b_\ell - a_\ell}{2} (y^{\ell,+}_0 + y^{\ell,-}_0) 
& \mbox{s.to}\quad & \vy^{\ell,+} - \vy^{\ell,-} = \vJ^\ell \vW^ \ell  \vp\\
\nonumber
& & \mbox{and}\quad & \Toep_\infty(\vy^{\ell,\pm}) \succeq {\bf 0}. 
\end{align}
This expression is not tractable due to the infinite dimensionality of the optimization variables and constraints,
but truncating them to a level $d$ leads to a tractable expression --- the above-mentioned ersatz norm.

\bprop
\label{ThmErs}
For each $d \ge N$, the expression \vspace{-1mm}
\begin{align}
\label{DefErs}
\sslash P \sslash_d & := 
\min_{\vy^{1,\pm},\ldots,\vy^{L,\pm} \in \bR^{d+1}}
\sum_{\ell=1}^L \f{b_\ell - a_\ell}{2} (y^{\ell,+}_0 + y^{\ell,-}_0) 
& \mbox{s.to}\quad & \vy^{\ell,+} - \vy^{\ell,-} = \vJ^{\ell,d} \vW^\ell \vp\\
\nonumber
& & \mbox{and}\quad & \Toep_d(\vy^{\ell,\pm}) \succeq {\bf 0}
\end{align}
defines a norm on the space of polynomials of degree at most $N$.
Moreover, one has
\be
\label{Monot}
\cdots \le \sslash P \sslash_d \le \sslash P \sslash_{d+1} \le \cdots \le \left\| \f{P}{w} \right\|_{L_1(K)}
\quad \mbox{and} \quad
\lim_{d \to \infty} \sslash P \sslash_d = \left\| \f{P}{w} \right\|_{L_1(K)}.
\ee
\eprop

\bpf
To justify that the expression in \eqref{DefErs} defines a norm,
we concentrate on the property $[\sslash P \sslash_d = 0] \imp [P=0]$,
as the other two norm  properties are fairly clear.
So, assuming that $\sslash P \sslash_d = 0$,
there exist $\vy^{1,\pm},\ldots,\vy^{L,\pm} \in \bR^{d+1}$ such that
\be
\label{=0}
\sum_{\ell=1}^L \f{b_\ell - a_\ell}{2} (y^{\ell,+}_0 + y^{\ell,-}_0)  = 0,
\ee
as well as, for all $\ell = 1:L$,
\be
\label{After=0}
\vy^{\ell,+} - \vy^{\ell,-} = \vJ^{\ell,d} \vW^ \ell \vp
\qquad \mbox{and} \qquad
\Toep_d(\vy^{\ell,\pm}) \succeq {\bf 0}.
\ee
The semidefiniteness of the Toeplitz matrices implies that
\be
 |y^{\ell,\pm}_k| \le y^{\ell,\pm}_0 
\qquad
\mbox{for all }k=0:d,
\ee
which, in view of \eqref{=0}, yields $\vy^{\ell,\pm} = {\bf 0}$.
By the invertibility of the matrices $\vW^\ell $ and the injectivity of the matrices $\vJ^{\ell,d}$ (easy to check from \eqref{EntriesJ}),
 we derive that $\vp = {\bf 0}$,
and in turn that $P=0$, as desired.

Let us turn to the justification of \eqref{Monot}.
The chain of inequalities translates the fact that the successive minimizations impose more and more constraints,
hence produce larger and larger minima.
\rev{It remains to prove} that the limit of the sequence $(\sslash P \sslash_d)_{d \ge N}$ equals $\|P/w\|_{L_1(K)}$ (the limit exists, because the sequence is nondecreasing and bounded above).
For each $d \ge N$,
as was done in \eqref{=0} and \eqref{After=0},
we consider minimizers of the problem \eqref{L1Reform2}
---
they belong to $\bR^{d+1}$ but we pad them with zeros to create infinite sequences $\vy^{1,\pm,d},\ldots,\vy^{L,\pm,d}$ satisfying
\be
\label{=0bis}
\sum_{\ell=1}^L \f{b_\ell - a_\ell}{2} (y^{\ell,+,d}_0 + y^{\ell,-,d}_0)  = \sslash P \sslash_d,
\ee
as well as, for all $\ell = 1:L$,
\be
\label{After=0bis}
\vy^{\ell,+,d} - \vy^{\ell,-,d} = \vJ^\ell \vW^ \ell \vp
\qquad \mbox{and} \qquad
\Toep_\infty(\vy^{\ell,\pm,d}) \succeq {\bf 0}.
\ee
The semidefiniteness of the Toeplitz matrices,
together with \eqref{=0bis},
implies that, for all $k \ge 0$,
\be
|y_k^{\ell,\pm,d}| \le y_0^{\ell,\pm,d} \le \f{2}{b_\ell - a_\ell} \sslash P \sslash_d
\le \f{2}{b_\ell - a_\ell}  \left\| \f{P}{w} \right\|_{L_1(K)}.
\ee
In other words,
each sequence  $(\vy^{\ell,\pm,d})_{d \ge N}$,
with entries in the sequence space $\ell_\infty$, is bounded. 
The sequential compactness Banach--Alaoglu theorem guarantees the existence of convergent subsequences in the weak-star topology.
With \rev{$(\vy^{\ell,\pm,d_m})_{m \ge 0}$ denoting these subsequences and}
$\vy^{\ell,\pm} \in \ell_\infty$ denoting their limits,
\rev{the weak-star convergence implies that}
\be
y^{\ell,\pm,d_m}_k \underset{m \to \infty}{\tto} y^{\ell,\pm}_k
\qquad \mbox{for all }k \ge 0.
\ee
Writing \eqref{After=0bis} for $d = d_m$ and passing to the limit reveals that the sequences $\vy^{1,\pm},\ldots,\vy^{L,\pm}$ are feasible for the problem \eqref{L1Reform2}.
Hence,
\begin{align}
\left\| \f{P}{w} \right\|_{L_1} &\le
\sum_{\ell = 1}^L \f{b_\ell - a_\ell}{2} (y^{\ell,+}_0 + y^{\ell,-}_0)
= \lim_{m \to \infty} \sum_{\ell = 1}^L \f{b_\ell - a_\ell}{2}  (y^{\ell,+,d_m}_0 + y^{\ell,-,d_m}_0)\\
\nonumber
& = \lim_{m \to \infty} \sslash P \sslash_{d_m}
= \lim_{d \to \infty} \sslash P \sslash_d,
\end{align}
\rev{where the last equality relied on the fact that the nondecreasing and bounded sequence $(\sslash P \sslash_d)_{d \ge N}$ is convergent.}
This concludes the justification of \eqref{Monot}.
\epf

Given $d \ge N$, let us now consider ersatz $N$th Chebyshev polynomial\rev{s} of the second kind \rev{for} $K$ (a priori not guaranteed to be unique) defined by
\be
\label{DefV}
\cV_{N,d}^K \rev{\in}  \underset{P(x) = x^N + \cdots }{\argmin} \sslash P \sslash_d.
\ee
It is possible to compute such a polynomial by solving the following semidefinite program:
\begin{align}
\label{MainSDP}
& \minimize{\substack{p_0, p_1, \ldots, p_N \in \bR \\ \vy^{1,\pm},\ldots,\vy^{L,\pm} \in \bR^{d+1}}}
\sum_{\ell=1}^L \f{b_\ell - a_\ell}{2} (y^{\ell,+}_0 + y^{\ell,-}_0) 
& \mbox{s.to}\quad & p_N = \f{1}{2^{N-1}}, \quad \vy^{\ell,+} - \vy^{\ell,-} = \vJ^{\ell,d} \vW^\ell \vp\\
\nonumber
& & \mbox{and}\quad & \Toep_d(\vy^{\ell,\pm}) \succeq {\bf 0}.
\end{align}
The qualitative result below ensures that, as $d$ increases,
the ersatz Chebyshev polynomial\rev{s} $\cV_{N,d}^K$ \rev{approach} genuine Chebyshev polynomial\rev{s} $\cU_{N}^K$,
which are themselves obtained by solving the following (unpractical) semidefinite program:
\begin{align}
\label{InfSDPforU}
& \minimize{\substack{p_0, p_1, \ldots, p_N \in \bR \\ \vy^{1,\pm},\ldots,\vy^{L,\pm} \in \bR^\bN}}
\sum_{\ell=1}^L \f{b_\ell - a_\ell}{2} (y^{\ell,+}_0 + y^{\ell,-}_0) 
& \mbox{s.to}\quad & p_N = \f{1}{2^{N-1}}, \quad \vy^{\ell,+} - \vy^{\ell,-} = \vJ^{\ell} \vW^\ell \vp\\
\nonumber
& & \mbox{and}\quad & \Toep_\infty(\vy^{\ell,\pm}) \succeq {\bf 0}.
\end{align}

\bthm
\label{ThmQual}
Any sequence $(\cV_{N,d}^K)_{d \ge N}$ of minimizers of \eqref{DefV} 
\rev{admits a subsequence converging
(with respect to any of the equivalent norms on the space of polynomials of degree at most $N$)
 to a minimizer $ \cU_N^K$ of \eqref{DefU}.
Moreover, if \eqref{DefU} has a unique minimizer $\cU_N^K$,
then the whole sequence  $(\cV_{N,d}^K)_{d \ge N}$ converges to $\cU_N^K$,
}
i.e.,
\be
\label{CV2}
\cV_{N,d}^K \underset{d \to \infty}{\tto} \cU_N^K.
\ee
\ethm

\bpf
We first  prove that the minima of \eqref{DefV} converge monotonically to the minimum of \eqref{DefU},
i.e.,
\be
\label{CV1}
\cdots \le \sslash \cV_{N,d}^K \sslash_d \le \sslash \cV_{N,d+1}^K \sslash_{d+1} \le \cdots \le \left\| \f{\cU_{N}^K}{w} \right\|_{L_1(K)}
\quad \mbox{and} \quad
\lim_{d \to \infty} \sslash \cV_{N,d}^K \sslash_d = \left\| \f{\cU_N^K}{w} \right\|_{L_1(K)}.
\ee
The argument is quite similar to the proof of \eqref{Monot} in Proposition \ref{ThmErs}.
The chain of \rev{inequalities} holds because more and more constraints are imposed.
Next, considering coefficients $p_0^d,p_1^d,\ldots,p_N^d$ and 
 infinite sequences $\vy^{1,\pm,d},\ldots,\vy^{L,\pm,d}$ satisfying
\be
\label{=0ter}
\sum_{\ell=1}^L \f{b_\ell - a_\ell}{2} (y^{\ell,+,d}_0 + y^{\ell,-,d}_0)  = \sslash \cV_{N,d}^K \sslash_d,
\ee
as well as $p^d_N = 1/2^{N-1}$ and, for all $\ell = 1:L$,
\be
\vy^{\ell,+,d} - \vy^{\ell,-,d} = \vJ^\ell \vW^ \ell \vp^d
\qquad \mbox{and} \qquad
\Toep_d(\vy^{\ell,\pm,d}) \succeq {\bf 0},
\ee
the semidefiniteness of the Toeplitz matrices, together with \eqref{=0ter}, still implies that the sequences $(\vy^{\ell,+,d})_{d \ge N}$ admit convergent subsequences in the weak-star topology,
so we can write
\be
y^{\ell,\pm,d_m}_k \underset{m \to \infty}{\tto} y^{\ell,\pm}_k
\qquad \mbox{for all }k \ge 0.
\ee
We note that
\be
\vp^{d_m} = (\vJ^{\ell,N} \vW^\ell)^{-1} (\vy^{\ell,+,d_m}_{\{0,\ldots,N\}} - \vy^{\ell,-,d_m}_{\{0,\ldots,N\}})
\underset{m \to \infty}{\tto} 
(\vJ^{\ell,N} \vW^\ell)^{-1} (\vy^{\ell,+}_{\{0,\ldots,N\}} - \vy^{\ell,-}_{\{0,\ldots,N\}}) =: \vp.
\ee
It is easy to see that the coefficients $p_0,p_1,\ldots,p_N \in \bR$ thus defined, together with the sequences $\vy^{1,\pm},\ldots,\vy^{L,\pm} \in \bR^\bN$, are feasible for the problem \eqref{InfSDPforU},
which implies that
\begin{align}
\left\| \f{\cU_N^K}{w} \right\|_{L_1} &\le
\sum_{\ell = 1}^L \f{b_\ell - a_\ell}{2} (y^{\ell,+}_0 + y^{\ell,-}_0)
= \lim_{m \to \infty} \sum_{\ell = 1}^L \f{b_\ell - a_\ell}{2}  (y^{\ell,+,d_m}_0 + y^{\ell,-,d_m}_0)\\
\nonumber
& = \lim_{m \to \infty} \sslash \cV_{N,d_m}^K \sslash_{d_m}
= \lim_{d \to \infty} \sslash \cV_{N,d}^K \sslash_d.
\end{align}
This concludes the justification of \eqref{CV1}.

\rev{Let us now prove that the sequence $(\cV_{N,d}^K)_{d \ge N}$ admits a subsequence converging to a minimizer $ \cU_N^K$ of \eqref{DefU}.
This sequence is bounded (with respect to any of the equivalent norms, e.g. $\sslash \cdot \sslash_{N}$):
indeed, 
as a consequence of \eqref{Monot} and \eqref{CV1}, we have
$\sslash \cV_{N,d}^K \sslash_N \le \sslash \cV_{N,d}^K \sslash_d \le \| \cU_N^K / w \|_{L_1(K)}$.
Therefore, there is a subsequence $(\cV_{N,d_m}^K)_{m \ge 0}$
converging to some monic polynomial $\cV_N^K$.
Let us assume that $\cV_N^K$ is not one of the minimizers $\cU_N^K$ of \eqref{DefU},
i.e., that $\|\cU_N^K/w\|_{L_1(K)} < \|\cV_N^K/w\|_{L_1(K)}$.}
In view of \eqref{Monot}, we can choose $d$ large enough so that
\be
\label{LastPf1}
\left\| \cV_N^K/w \right\|_{L_1(K)}
< \sslash \cV_N^K \sslash_d + \eps,
\qquad \mbox{where} \quad
\eps := \left\| \cV_N^K/w \right\|_{L_1(K)} - \left\| \cU_N^K/w \right\|_{L_1(K)} > 0.
\ee
Let us observe that,
\rev{with $d$ being fixed and}
by virtue of \eqref{Monot} and \eqref{CV1}, 
\be
\label{LastPf2}
\sslash \cV_N^K \sslash_d 
= \lim_{m \to \infty} \sslash \cV_{N,d_m}^K \sslash_d
\le \lim_{m \to \infty} \sslash \cV_{N,d_m}^K \sslash_{d_m}
= \left\| \cU_N^K/w \right\|_{L_1(K)} .
\ee
Combining \eqref{LastPf1} and \eqref{LastPf2} yields
\be
\left\| \cV_N^K/w \right\|_{L_1(K)}
< \left\| \cU_N^K/w \right\|_{L_1(K)} + \eps
= \left\| \cV_N^K/w \right\|_{L_1(K)},
\ee
which is of course a contradiction.
\rev{This implies that $\cV_N^K$ is a minimizer of \eqref{DefU}, as expected.}\newpage

\rev{Finally, in case \eqref{DefU} has a unique minimizer $\cU_N^K$, 
we can establish \eqref{CV2} by contradiction.
Namely, if the sequence $(\cV_{N,d}^K)_{d \ge N}$ did not converge to $\cU_{N}^K$,
then we could construct a subsequence $(\cV_{N,d_m}^K)_{m \ge 0}$
converging to some monic polynomial $\cV_N^K \not= \cU_N^K$.
Repeating the above arguments would imply that $\cV_N^K$ is a minimizer of \eqref{DefU}, i.e., $\cV_N^K = \cU_N^K$,
providing the required contradiction.}
\epf

Theorem \ref{ThmQual} does not indicate how to \rev{choose $d$ a priori} in order to reach a prescribed accuracy for the distance between $\cV_{N,d}^K$ and $\cU_N^K$.
However, for a given $d$,
we can assess a posteriori the distance between the ersatz minimum $\sslash \cV_{N,d}^K \sslash_d$ and the genuine minimum $\|\cU_N^K/w\|_{L_1(K)}$.
Indeed, on the one hand, the semidefinite program \eqref{MainSDP} produces $\sslash \cV_{N,d}^K \sslash_d$ while outputting $\cV_{N,d}^K$;
on the other hand, the weighted $L_1$-norm $\|\cV_{N,d}^K/w\|_{L_1(K)}$ can be computed once $\cV_{N,d}^K$ has been output.
These two facts provide lower and upper bounds for the unknown value $\|\cU_N^K/w\|_{L_1(K)}$,
as stated by the quantitative result below.

\bprop
For any $d \ge N$, one has
\be
\label{Assess}
\sslash \cV_{N,d}^K \sslash_d
 \le  \left\| \cU_N^K / w \right\|_{L_1(K)} 
 \le \left\| \cV_{N,d}^K / w \right\|_{L_1(K)},
\ee
hence the weighted $L_1$-norm of $\cU_N^K$ on $K$ is approximated with a computable relative error of
\be 
\delta_{N,d}^K := 1-\f{\sslash \cV_{N,d}^K \sslash_d}{\| \cV_{N,d}^K/w \|_{L_1(K)}} \ge 0.
\ee
\eprop

\bpf
By the definition \eqref{DefU} of the genuine Chebyshev polynomial of the second kind,
we have
\be
\|\cU_N^K/w\|_{L_1(K)} \le \|\cV_{N,d}^K/w\|_{L_1(K)},
\ee
and by the definition \eqref{DefV} of the ersatz Chebyshev polynomial of the second kind, together with~\eqref{Monot}, we have
\be
\sslash \cV_{N,d}^K \sslash_d
\le \sslash \cU_N^K \sslash_d
\le \|\cU_N^K/w\|_{L_1(K)}.
\ee
This establishes the bounds announced in \eqref{Assess}.
We also notice that the relative error satisfies
\be
\label{dTo0}
\delta_{N,d}^K = \f{\|\cV_{N,d}^K/w\|_{L_1(K)} - \sslash \cV_{N,d}^K \sslash_d}{\|\cV_{N,d}^K/w\|_{L_1(K)}}
 \underset{d \to \infty}{\tto} 0,
\ee
since, \rev{according to} \eqref{CV2} and \eqref{CV1},
both $\|\cV_{N,d}^K/w\|_{L_1(K)}$ and $\sslash \cV_{N,d}^K \sslash_d$ converge to 
$\|\cU_N^K/w\|_{L_1(K)}$
\rev{in case of uniqueness of $\cU_N^K$.
In case of nonunuqueness, \eqref{dTo0} remains true at least for a subsequence.}
\epf

\begin{figure}[hbtp]
\center
\subfigure[
\rev{$K=K_1$: $\delta_{N,d}^K \approx 7 \cdot 10^{-4}$}
]{
\includegraphics[width=0.48\textwidth,height=0.33\textwidth]{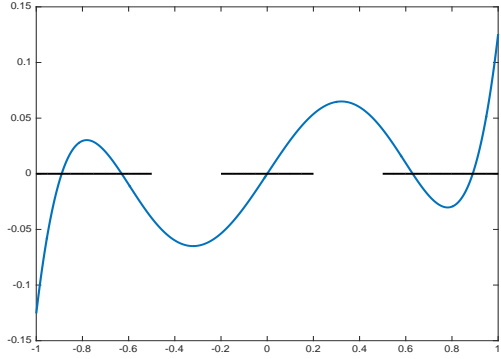}
}
\subfigure[
\rev{$K=K_1$, weighted: $\delta_{N,d}^K \approx 6 \cdot 10^{-4}$}
]{
\includegraphics[width=0.48\textwidth,height=0.33\textwidth]{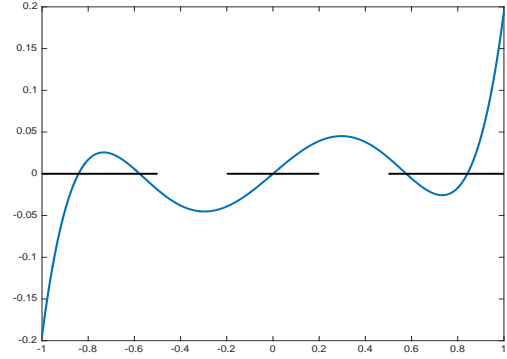}
}

\subfigure[
\rev{$K=K_2$: $\delta_{N,d}^K \approx 8 \cdot 10^{-4}$}
]{
\includegraphics[width=0.48\textwidth,height=0.33\textwidth]{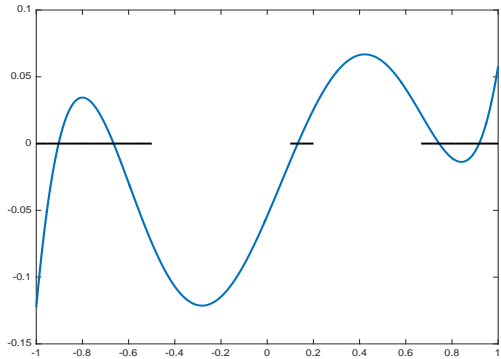}
}
\subfigure[
\rev{$K=K_2$, weighted: $\delta_{N,d}^K \approx 8 \cdot 10^{-4}$]}
]{
\includegraphics[width=0.48\textwidth,height=0.33\textwidth]{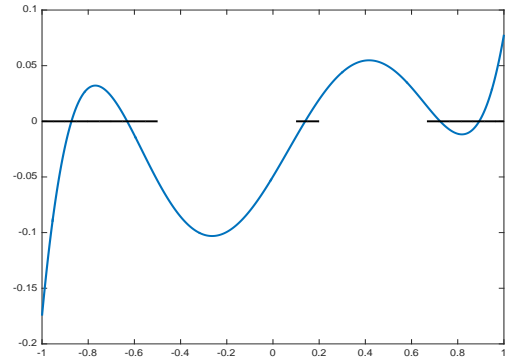}
}
\caption{
\rev{Ersatz $5$th Chebyshev polynomials of the second kind for $K_1 = \big[-1,-\f{1}{2}\big] \cup \big[-\f{1}{5},\f{1}{5} \big] \cup \big[\f{1}{2},1 \big]$ and for $K_2 = \big[ -1,-\f{1}{2} \big] \cup \big[\f{1}{10},\f{1}{5}\big] \cup \big[\f{2}{3},1\big]$:
the first column corresponds to the unweighted case,
while weighted ersatz Chebyshev polynomials with weight $w(x) = (1+x^2)/(2-x^2)$ are shown in the second column.}
}
\label{Fig2}
\end{figure}

Figure \ref{Fig2} shows ersatz Chebyshev polynomials of the second kind computed on the same examples as in Figure \ref{Fig1}.
Notice that no `restricted' ersatz Chebyshev polynomials of the second kind are displayed.
This is because our experiments suggested that the polynomials $\cV_{N,d}^K$ had \rev{simple roots all lying inside $K$.}
The corresponding statement for the polynomials $\cU_N^K$\rev{, in case of uniqueness,} can in fact be justified theoretically \rev{by the following observation}.

\bprop
\rev{Let $\cU_N^K$ be a weighted Chebyshev polynomial of the second kind for a finite union of closed intervals $K \inc [-1,1]$.
This polynomial is the unique minimizer of \eqref{DefU} if and only if it has $N$ simple roots all lying inside $K$.}
\eprop

\bpf
As a minimizer of \eqref{DefU},
\rev{a} Chebyshev polynomial of the second kind \rev{for} $K$ is characterized (see e.g. \cite[p. 84, Theorem 10.4]{CA}) by the condition
\be
\label{Cond}
\int_{K} \f{\rev{\sgn( \cU_N^K(x) ) P(x)}}{w(x)} \, dx = 0
\qquad \qquad \mbox{for all polynomials $P$ of degree less than $N$}.
\ee
This implies that $\cU_{N}^K$ has $N$ roots in \rev{$(-1,1)$},
as \eqref{Cond} \rev{would} not hold \rev{for} $P(x) = (x-\xi_1) \cdots (x-\xi_n)$
if $\cU_N^K$ had $n<N$ roots $\xi_1,\ldots,\xi_n$ in \rev{$(-1,1)$}.
\rev{Moreover, if one of the roots was repeated, 
we would have $\cU_N^K(x) = (x-\xi)^2 P(x)$ for some polynomial $P$ of degree $N-2$,
but then  \eqref{Cond} would not hold for this $P$ either.
Thus, the polynomial $\cU_N^K$ can be written, with distinct $\xi_1, \ldots, \xi_N \in (-1,1)$, as
\be
\cU_N^K(x) = (x-\xi_1) \cdots  (x-\xi_i) \cdots (x-\xi_N).
\ee
Assume that $\cU_N^K$ is the unique Chebyshev polynomial of the second kind for $K$.
If one of the $\xi_i$'s does not lie inside $K$,
i.e., if $\xi_i$ belongs to one of the gaps $[b_\ell,a_{\ell+1}]$,
then we can perturb $\xi_i$ to $\wt{\xi}_i$ while keeping it in $[b_\ell,a_{\ell+1}]$.
Hence, the perturbed monic polynomial 
\be
\wt{\cU}_N^K(x) =  (x-\xi_1) \cdots  (x-\wt{\xi}_i) \cdots (x-\xi_N).
\ee
still satisfies $\sgn(\wt{\cU}_{N}^K(x)) = \sgn( \cU_N^K(x))$ for all $x \in K$.
The condition \eqref{Cond} is then fulfilled by $\wt{\cU}_{N}^K$, too, 
so this monic polynomial is another minimizer of \eqref{DefU}, which is impossible.
We have therefore proved that the $N$ simple roots of $\cU_N^K$ all lie inside $K$.

Conversely, assume that $\cU_N^K$ has $N$ simple roots all lying inside $K$
and let us prove that $\cU_N^K$ is the unique minimizer of \eqref{DefU}.
Consider a monic polynomial $\wt{\cU}_N^K $ 
with $\|\wt{\cU}_N^K / w\|_{L_1(K)} = \|\cU_N^K / w\|_{L_1(K)}$.
In view of \eqref{Cond}, we notice that
\be
\int_{K} \f{\sgn( \cU_N^K(x) ) (\cU_N^K(x) - \wt{\cU}_N^K(x))}{w(x)} \, dx = 0.
\ee
From here,
it follows that
\begin{align}
\bigg\| \f{\cU_N^K}{w} \bigg\|_{L_1(K)}
& = \int_K \f{|\cU_N^K(x)|}{w(x)} \, dx
= \int_{K} \f{\sgn( \cU_N^K(x) ) \cU_N^K(x) }{w(x)} \, dx\\
\nonumber
& = \int_{K} \f{\sgn( \cU_N^K(x) ) \wt{\cU}_N^K(x)}{w(x)} \, dx
\le \int_K \f{|\wt{\cU}_N^K(x)|}{w(x)} \, dx
 = \bigg\| \f{\wt{\cU}_N^K}{w} \bigg\|_{L_1(K)}.
\end{align}
 The first and the last terms being equal,
we must have equality all the way through, 
which means that $\sgn( \wt{\cU}_N^K(x)) = \sgn( \cU_N^K(x))$ for all $x \in K$.
Given that the polynomial $\cU_N^K$ vanishes at distinct points $\xi_1,\ldots,\xi_N$ inside  $K$,
the polynomial $\wt{\cU}_N^K$ must also vanish at $\xi_1,\ldots,\xi_N$,
and since both polynomials are monic,
we must have $\wt{\cU}_N^K = \cU_N^K$,
proving the uniqueness.}
\epf

\end{document}